\documentclass[letterpaper,onecolumn,11pt]{ieeeconf}
\IEEEoverridecommandlockouts 
\usepackage{amsmath}
\usepackage{amssymb}
\usepackage{ntheorem}
\usepackage{color}
\usepackage{bbold}
\usepackage{epsfig}
\usepackage{graphicx}
\usepackage{url}

\newtheorem{remark}{Remark}

\newcommand{\tr}{\intercal} 
\newcommand{\defeq}{\overset{\mathrm{def}}{=}}

\title{\bf Distributed MPC with ALADIN---A Tutorial}

\author{Boris Houska, Jiahe Shi
\thanks{Boris Houska and Jiahe Shi are with the School of Information Science and Technology, ShanghaiTech University, China;
{\tt borish, shijh@shanghaitech.edu.cn}.}%
}

\begin{document}

\maketitle

\begin{abstract}
This paper consists of a tutorial on the Augmented Lagrangian based Alternating Direction Inexact Newton method (ALADIN) and its application to distributed model predictive control (MPC). The focus is---for simplicity of presentation---on convex quadratic programming (QP) formulations of MPC. It is explained how ALADIN can be used to synthesize sparse QP solvers for large-scale linear-quadratic optimal control by combining ideas from augmented Lagrangian methods, sequential quadratic programming, as well as barrier or interior point methods. The highlight of this tutorial is a real-time ALADIN variant that can be implemented with a few lines of code yet arriving at a sparse QP solver that can compete with mature open-source and commercial QP solvers in terms of both run-time as well as numerical accuracy. It is discussed why this observation could have far reaching consequences on the future of algorithm and software development in the field of large-scale optimization and MPC.  
\end{abstract}

\section{Introduction}
The success of MPC~\cite{Rawlings2017} in academic and industrial applications~\cite{Qin2003} relies on high-performance real-time optimization algorithms. Over the last decades, such real-time optimization algorithms and software have been developed for  small- to medium-scale optimization problems, which can solve these problems online within the milli- and microsecond range~\cite{Houska2011,Mattingley2012}. Numerical inaccuracies or the run-time of such solvers are usually only a problem if one attempts to implement MPC for a system with a very large number of states or controls, or if one formulates MPC problems for nonlinear models.

The goal of the present tutorial paper is not only to review the state-of-the-art but also to discuss recent achievements concerning the development of real-time distributed optimization methods for MPC based on ALADIN~\cite{Houska2016,Houska2021}. Here, our focus is, for simplicity of presentation, on linear systems with quadratic stage cost such that the MPC problem can be formulated as a convex QP. Before being able to understand why this focus on such a basic MPC setting is of interest---despite the fact that QP solvers for MPC have been developed over the past $70$ years and despite the fact that many generic and tailored QP solvers have reached a high level of maturity---we need to briefly review existing solution methods for QPs. Namely, there are three big classes of existing methods: active set methods~\cite{Wolfe1959}, as implemented in the software packages \texttt{qpOASES}~\cite{Ferreau2014}, \texttt{MOSEK}~\cite{Mosek2022} or \texttt{GUROBI}~\cite{Gurobi2022}, interior point methods~\cite{Nesterov1994}, as implemented in~\texttt{CVXGEN}~\cite{Mattingley2012} and \texttt{OOQP}~\cite{Gertz2003}, as well as first order \mbox{methods~\cite{Frank1956}---nowadays} often based on the alternating direction method of multipliers \mbox{(ADMM)~\cite{Boyd2011}---as,} for example, implemented in the open-source conic solver \texttt{SCS}~\cite{Donoghue2016} as well as the sparse QP solver \texttt{OSQP}~\cite{Stellato2020}. Notice that most of the latter references come along with more complete overviews of the history of QP solver development to which we refer at this point.

In contrast to the above reviewed convex QP solvers, ALADIN has originally been developed for solving large-scale non-convex optimization problems. For instance, the original ALADIN variant~\cite{Houska2016} combines ideas from the field of sequential quadratic programming (SQP)~\cite{Nocedal2006} and augmented Lagrangian methods~\cite{Hamdi2011} in order to compose a second order method for distributed non-convex optimization. Nevertheless, certain variants of ALADIN can also be used as first order methods for non-differentiable convex optimization problems. In such a setting global convergence of ALADIN can be established~\cite{Houska2021}.

\subsection*{Highlights}
As the present paper consists of a tutorial on how to apply existing variants of ALADIN to develop sparse and distributed QP solvers for MPC, no new theoretical results are presented. Nevertheless, after a short introduction to MPC in Section~\ref{sec::mpc}, the following two new aspects of ALADIN can be considered as highlights of this tutorial.

\begin{enumerate}

\item Section~\ref{sec::ALADIN} introduces a complete and very practical variant of ALADIN for solving sparse and distributed QPs. This variant is based on a synthesis of active set methods, interior point methods, and first order methods, which leads to a new high-performance sparse QP solver that can be implemented with just a few lines of code.

\item Section~\ref{sec::RealTime} reviews ideas from~\cite{Jiang2021} on how to develop real-time distributed variants of the presented sparse QP solver that are tailored for MPC. An application of this solver to large-scale benchmark MPC problems illustrates the competitiveness of ALADIN compared to other existing QP solvers.

\end{enumerate}

As much as this tutorial explains how to synthesize a ``simple'' variant of ALADIN that can compete with existing sparse QP solvers, our goal is not to develop yet another QP solver software package. On the contrary, the ultimate goal of this line of research on ALADIN is of a completely different nature: ALADIN has originally been designed for solving large-scale non-convex optimization and nonlinear MPC problems. At the current status of research, very early-stage software packages based on ALADIN have appeared~\cite{Engelmann2022}. However, the step from such early-stage implementations to a high-performance large-scale non-convex problem solver will require much research effort and time investment. Thus, for numerical software developers who wish to work on such large-scale optimization software, it will be important to assess first whether investing time into ALADIN based solvers has the potential to advance the state-of-the-art. The present tutorial intends to help with the porgress to come to such an assessment by pointing out that ALADIN works well for solving convex QPs. This can be interpreted as one possible indicator that further research on the implementation of a non-convex optimization problem solver could indeed be fruitful. Therefore, Section~\ref{sec::conclusion} will not only summarize the highlights of this tutorial, but also elaborate on how the numerical observations from this article might impact the future of high-performance numerical optimization solver development for convex and non-convex programming as well as distributed MPC.

\subsection{Notation}
Let $\mathbb R^n$ denote the $n$-dimensional real vector space and
\[
\Vert x \Vert_Q^2 \ \defeq \ x^\tr Q x
\]
a weighted Euclidean norm with $Q \in \mathbb R^{n \times n}$ denoting a positive definite matrix. We occasionally use the notation
\[
\left\| x \right\|_{Q_1,Q_2}^2 \ \defeq \Vert x_1 \Vert_{Q_1}^2 + \Vert x_2 \Vert_{Q_2}^2 ,
\]
where $x = \left[x_1^\tr,x_2^\tr \right]^\tr$ can be a block vector of any dimension and $Q_1$ and $Q_2$ denote positive matrices of the corresponding dimensions. The symbols $\mathbb{1}$ and $\mathbf{1}$ are used to denote, respectively, the unit matrix and a vector whose components are all equal to one---assuming that it is clear from the context what their dimensions are.

\section{Model Predictive Control}
\label{sec::mpc}

This section reviews linear-quadratic MPC controllers.

\subsection{Linear-Quadratic MPC}
This paper concerns MPC problems of the form 
\begin{align}
J_N(x_0) \ \defeq \ \min_{x,u} \ & \sum_{k=0}^{N-1} \left\{ \Vert x_k \Vert_Q^2 + \Vert u_k \Vert_R^2 \right\} + \Vert x_N \Vert_P^2 \notag \\[0.1cm]
\label{eq::MPC}
\mathrm{s.t.} \ & \left\{
\begin{array}{l}
\forall k \in \{ 0,1,\ldots, N-1 \}, \\[0.1cm]
x_{k+1} = A x_k + B u_k \\[0.1cm]
c \leq C x_k + D u_k \leq d \; .
\end{array}
\right.
\end{align}
Here, $x_k \in \mathbb R^{n_x}$ denotes the state and $u_k \in \mathbb R^{n_u}$ the control input at time $k$. The corresponding MPC feedback law,
\[
\mu(x_0) = u_0^\star(x_0),
\]
corresponds to the first element, $u_0^\star(x_0)$, of the minimizing control input sequence of~\eqref{eq::MPC}, which depends on the state measurement $x_0$. Throughout this paper the system matrices $A$ and $B$ as well as the joint state- and control constraint matrices $C$ and $D$ are assumed to be given. Moreover, for simplicity of presentation, we assume that the matrices $Q$, $R$, and $P$ are positive definite, although several of the considerations below can be generalized for positive semi-definite weights, too. As reviewed in the introduction, there exist many numerical methods for solving~\eqref{eq::MPC}. Nevertheless, numerical challenges can arise if we have a system with a very large number of states, $n_x \gg 1$, while $A$ and $Q$ (and sometimes also $B$ and $R$) are sparse matrices. Additionally, in some applications, for instance, if $n_u \ll n_x$, one  needs large prediction horizons $N$ in order to achieve a satisfying control performance.

\begin{remark}
A historical overview, recent developments, and a more complete discussion on how to write economic and distributed MPC problems in the form of~\eqref{eq::MPC} can be found in the book~\cite{Rawlings2017} and the recent overview article~\cite{Muller2017}.
\end{remark}

\subsection{Recursive Feasibility and Stability}
For simplicity of presentation, we assume that $(A,B)$ is asymptotically stabilizable such that we can compute the positive definite matrix $P$ by solving the Riccati equation
\begin{align}
\label{eq::riccati}
P =  A^\tr P A + Q - A^\tr P B (R + B^\tr PB)^{-1}B^\tr P A.
\end{align}
We additionally assume that $c < 0 < d$ such that all constraints are strictly feasible. Thus, if we choose $N$ sufficiently large, we have $J_N(x_0) = J_\infty(x_0)$ for all $x_0$ that are in the domain of the infinite horizon cost $J_\infty$. Here, $J_\infty$ is a piecewise quadratic and positive definite function that satisfies the stationary Bellman equation 
\begin{eqnarray}
\label{eq::Bellman}
\begin{array}{rcl}
J_\infty(x) = &\min_{u}& \ \Vert x \Vert_Q^2 +  \Vert u \Vert_R^2 + J_\infty(Ax+Bu) \\[0.16cm]
&\mathrm{s.t.}& \ c \leq Cx+Du \leq d \; .
\end{array}
\end{eqnarray}
Notice that under these assumption the MPC controller is recursively feasible and $\mu$ stabilizes the system~\cite{Rawlings2017}. However, such stability and recursive feasibility statements only hold if~\eqref{eq::MPC} is solved exactly. For real-time MPC solvers one needs to impose additional requirements on the accuracy of the solver in order to ensure stability; see Section~\ref{sec::RealTime}.

\subsection{Tutorial Example}
\label{sec::tutorial}
Throughout this paper, we use a chained spring-mass-damper system as a tutorial. The position and velocity of the $i$-th wagon, $p_i$ and $v_i$, satisfy a recursion of the form
\begin{eqnarray}
p_i^+ &=& p_i + h v_i \notag \\
v_i^+ &=& v_i + h \left[ \frac{k_\mathrm{s}}{m} ( p_{i-1} - 2 p_{i} + p_{i+1} ) -  \frac{k_\mathrm{d}}{m} v_i + \frac{u_i}{m} \right] \, \notag
\end{eqnarray}
for all $i \in \{ 1, \ldots, \mathsf{n} \}$. Here, we formally define $p_{0} = 0$ as well as $p_{\mathsf{n}+1} = p_\mathsf{n}$, which means that the first wagon of the chain is attached to a wall while the last wagon is free. Moreover, $u_i$ denotes a piecewise constant force at the $i$-th wagon. For simplicity of presentation, we set the spring constant, the damping constant and mass to \mbox{$k_\mathrm{s} = k_\mathrm{d} = m = 1$}. Moreover, we use the Euler discretization parameter $h = 0.1$. The matrices and vectors
\[
C = \left( \hspace{-0.1cm}
\begin{array}{cc}
\mathbb 1 & 0 \\
0 & 0
\end{array} \hspace{-0.1cm}
\right), \ D = \left( \hspace{-0.1cm}
\begin{array}{cc}
0 & 0 \\
0 & \mathbb 1 
\end{array} \hspace{-0.1cm}
\right), \ \ \text{and} \ \ d = -c = \left( \hspace{-0.1cm}
\begin{array}{c}
5 \cdot \mathbf{1} \\
\mathbf{1}
\end{array} \hspace{-0.1cm}
\right)
\]
model simple state- and control constraints. Similarly, we set $Q = \mathbb{1}$ and $R = \mathbb{1}$ while $P$ is computed by solving the above mentioned algebraic Riccati equation.

\section{ALADIN}
\label{sec::ALADIN}

This section explains how ALADIN can be used as a generic sparse QP solver for solving~\eqref{eq::MPC}. 

\subsection{Distributed Quadratic Programming}

There are several ways to exploit the structure of~\eqref{eq::MPC}. Similar to the implementation of ADMM in the software package \texttt{OSQP}~\cite{Stellato2020}, we focus in this paper on a generic sparse structure exploitation scheme writing~\eqref{eq::MPC} in the form
\begin{align}
\label{eq::qp}
\begin{array}{cl}
\underset{y,z}{\min} & F(y) + G(z) \\[0.16cm]
\mathrm{s.t.} & E y = z \ \mid \lambda \; ,
\end{array}
\end{align}
where the auxiliary vector $y = [u_0^\tr,x_1^\tr,u_1^\tr,\ldots,x_N^\tr]^\tr$ collects all optimization variables of the MPC problem and $\lambda$ denotes the dual solution. Here, $F$ is a quadratic function,
\begin{eqnarray}
\notag
F(y) = \frac{1}{2} y^\tr \mathcal Q y &\defeq& \sum_{k=0}^{N-1} \left\{ \Vert x_k \Vert_Q^2 + \Vert u_k \Vert_R^2 \right\} + \Vert x_N \Vert_P^2 \; .
\end{eqnarray}
The sparse matrix $E$ is given by
\[
E \ \defeq \ \left(
\begin{array}{rrrrrrrrr}
B & -\mathbb{1} &  &  &  &  &  &  \\
D & 0 &  &  &  &  &  &  \\
    & A & B & -\mathbb{1} &  &   &  & \\
    & C & D & 0 &  &   &  &  \\
    &  &  &  & \ddots & &  \\
    &  &  &   &  & A & B & -\mathbb{1} \\
    &  &  &   &  & C & D & 0
\end{array}
\right),
\]
where the empty blocks are all equal to $0$. Notice that~\eqref{eq::MPC} and~\eqref{eq::qp} are equivalent, if we define the function $G$ as  
\begin{align}
\label{eq::G}
& G(z) \ \defeq \ \left\{
\begin{array}{ll}
0 & \text{if} \; \underline z
\leq z \leq \overline z \\[0.1cm]
\infty & \text{otherwise}
\end{array}
\, \right\} \\[0.1cm]
& \text{with} \quad \underline z \ \defeq \  [-(A x_0)^\tr,(c-Cx_0)^\tr, 0^\tr, c^\tr, \ldots, 0^\tr,c^\tr]^\tr \notag \\[0.16cm]
& \text{and} \quad \ \overline z \ \defeq \ [-(A x_0)^\tr,(d-Cx_0)^\tr, 0^\tr, d^\tr \ldots, 0^\tr,d^\tr]^\tr \; .
\notag
\end{align}
The parametric initial value $x_0$ enters via the vectors $\underline z$ and $\overline z$ while all matrices are constant.

\subsection{Distributed Optimization Algorithm}

The main idea of ALADIN is to start with an initial guess $(y,z)$ for the primal solution of~\eqref{eq::qp} as well as an initial guess $\lambda$ for its dual solution and repeat the following steps.

\begin{enumerate}
\addtolength{\itemsep}{2pt}

\item Choose positive definite matrices $H \succ 0$, $\Sigma \succ 0$, and $K \succ 0$ as well as a tuning parameter $\theta \in [0,1]$.

\item Solve the decoupled optimization problems
\begin{align}
\label{eq::dco1}
& \min_v \; F(v) + \lambda^\tr E^\tr v + \frac{1}{2} \Vert v-y \Vert_{\Sigma}^2 \\[0.1cm]
\label{eq::dco2}
\text{and} \quad & \min_w \; G(w) - \lambda^\tr w + \frac{1}{2} \Vert w-z \Vert_{K}^2
\end{align}
and denote the minimizers by $v$ and $w$.

\item Compute the gradient of $F$ and subgradient of $G$ as
\begin{align}
\label{eq::gradF}
\nabla F(v) &= \Sigma(y-v) - E^\tr \lambda  \\[0.1cm]
\label{eq::gradG}
\text{and} \quad \partial G(w) &= K(z-w) + \lambda
\end{align}
at the decoupled minimizers $v$ and $w$.

\item Solve the equality constrained consensus QP
\begin{align}
\underset{y^+,z^+}{\min} \ & \frac{1}{2} \left\|
\left( \hspace{-0.1cm}
\begin{array}{c}
y^+-v \\
z^+-w
\end{array} \hspace{-0.1cm}
\right)
\right\|_{H,K}^2 + \left( \hspace{-0.1cm}
\label{eq::consensus}
\begin{array}{c}
\nabla F(v) \\
\partial G(w)
\end{array}
\hspace{-0.1cm} \right)^\tr
\left( \hspace{-0.1cm}
\begin{array}{c}
y^+ \\
z^+
\end{array}
\hspace{-0.1cm} \right)
\notag \\[0.16cm]
\mathrm{s.t.} \ & E y^+ = z^+ \; \mid \; \lambda^+ \; ,
\end{align}
denote the primal minimizers by $y^+$ and $z^+$, and denote the dual solution by $\lambda^+$.

\item Update the variables $y \leftarrow y^+$, $z \leftarrow z^+$, as well as $\lambda \leftarrow \theta \lambda^+ + (1-\theta) \partial G(y)$; and go to Step 1.

\end{enumerate}

\noindent
Notice that, if we set $\theta = 1$, the above algorithm coincides with the derivative-free variant of ALADIN that has been analyzed in~\cite{Houska2021}. However, the corresponding global convergence proof can be generalized easily for any choice of $\theta \in [0,1]$. As such, the above algorithm converges for all convex QPs---even without requiring positive definiteness of the objective matrices.

\begin{remark}
The above algorithm is neither equivalent to SQP nor to ADMM, although the introduction of the decoupled augmented Lagrangian problems in Step~2) is inspired by ADMM, while Step~10) is inspired by SQP. In fact, if $F$ and $G$ would both be smooth, one could---in complete analogy to SQP methods---set $H$ and $K$ to the Hessian matrices of $F$ and $G$ in order to obtain a variant of ALADIN that has a locally quadratic convergence rate~\cite{Houska2016}. In our context, however, $G$ is a non-smooth  function. Thus, the choice of $K$ requires further discussion.
\end{remark}

\subsection{Hessian Matrix Updates}

In our context, $F$ is a smooth quadratic form. Consequently, we set the Hessian matrix approximation $H$ to
\[
H \ \defeq \mathcal Q + \epsilon_1 \cdot \mathbb{1} \; ,
\]
which coincides with the exact Hessian of $F$ apart from a small regularization term that can be adjusted by the tuning parameter $\epsilon_1 > 0$. The reason for introducing this tuning parameter is that, for $\epsilon_1 > 0$, it can be shown that the above outlined variant of ALADIN also converges for LPs, where we have $\mathcal Q = 0$, or other types of degenerate QPs, where $\mathcal Q$ might be highly ill-conditioned~\cite{Houska2021}.

Next, in order to be able to assign a ``Hessian approximation'' to the non-smooth function $G$, one option is to establish an analogy to traditional interior point methods for convex optimization. In order to elaborate on this idea, we introduce relaxed log-barrier functions of the form
\[
\varphi(r,t,\xi) \, = \, -\frac{1}{t} \cdot \log \left( r-\xi \right)
\]
with barrier parameter \mbox{$t > 0$} and relaxation parameter \mbox{$r \geq 0$}, which are defined on the domain $\xi < r$. An associated smooth approximation, $G(z) \approx \Phi(r,t,z)$, is then given by
\begin{eqnarray}
\Phi( r,t, z ) &\defeq& \sum_{i} \, \left( \; \varphi \left( r,t, \underline z_i-z_i \right) + \varphi \left( r,t, z_i-\overline z_i \right) \; \right) . \notag
\end{eqnarray}
This approximation becomes exact for increasing barrier parameters and vanishing relaxation, $t \to \infty$ and $r \to 0^+$. In the following implementation of ALADIN, we use the Hessian matrix approximation $K = \nabla_z^2 \Phi(r,t,z)$, which is motivated by the fact that $G \approx \Phi(r,t,\cdot)$. An explicit expression for $K$ can be found by using the formula
\[
\nabla_{\xi}^2 \, \varphi(r,t,\xi) \ = \ \frac{1}{t} \frac{1}{(r-\xi)^2} \ > \ 0
\]
when evaluating the second order derivative of $\Phi$. Notice that the log-barrier is here merely used for tuning the matrix $K$ by interpreting this matrix as a Hessian approximation. In contrast to actual interior point methods, however, the decoupled optimization problem~\eqref{eq::dco2} uses the exact function $G$ rather than its relaxed log-barrier approximation.

\subsection{Implementation Details}
Notice that almost all steps of the above outlined algorithm involve the solution of equality constrained convex QPs, which can be solved by using sparse Cholesky factorizations.  The only exception is the decoupled optimization problem~\eqref{eq::dco2}. Since the above outlined log-barrier based weight matrix generation scheme leads to  a diagonal $K$, this problem can, however, be solved by a projection onto the box $[\underline z,\overline z]$. More precisely, if $\Pi$ denotes the projection function,
\[
\forall i \in \{ 1, \ldots, n_E \}, \quad \Pi_i(\xi) = \left\{
\begin{array}{ll}
\underline z_i & \text{if} \quad \xi_i < \underline z_i \\[0.1cm]
\xi_i \quad & \text{if} \quad \underline z_i \leq \xi_i \leq \overline z_i \\[0.1cm]
\overline z_i & \text{if} \quad \overline z_i < \xi_i \; ,
\end{array}
\right.
\]
the explicit solution for the decoupled variable $v$ is given by $v = \Pi( z + K^{-1} \lambda )$. With this, we have all ingredients to setup a complete sparse QP solution algorithm, as summarized in Figure~\ref{fig::alg1}.
\begin{figure}
\hrule
\hrule
\vspace{0.3cm}
\begin{center}
\textbf{ALADIN as Sparse QP Solver}
\end{center}

\hrule

\smallskip

\begin{itemize}
\addtolength{\itemsep}{3pt}
\item \textbf{Default Regularization and Initialization:}
\begin{itemize}
\addtolength{\itemsep}{3pt}
\item Set $\epsilon_1 = 10^{-6}$ and $\epsilon_2 = 10^{-3}$.
\item If the user does not specify any customized initialization, set $y = 0$ and $z=\lambda=0$.
\end{itemize}

\item \textbf{Initial Hessian Approximations:}
\begin{itemize}
\addtolength{\itemsep}{3pt}
\item Set $H = \mathcal Q + \epsilon_1 \cdot \mathbb{1}$, $\Sigma = \mathcal Q + \epsilon_2 \cdot \mathbb{1}$, and $K = \mathbb{1}$.
\end{itemize}

\item \textbf{Main Loop:}

\item[] \textbf{For $\; i = 1:i_\mathrm{max} \;$ do:}

\begin{enumerate}

\item Set $\; \sigma \ \leftarrow \ E^\tr \lambda$.

\item Set $\; v \ \leftarrow \ [ \mathcal Q+\Sigma]^{-1}( \Sigma y - \sigma )$.

\item Set $\; w \ \leftarrow \ \Pi( z + K^{-1} \lambda )$.

\item If $\| w - z \|_\infty \leq \mathsf{TOL}$ and $\| \mathcal Q y + \sigma \|_\infty \leq \mathsf{TOL}$, \textbf{break}.

\item Set $\; h \ \leftarrow \ \Sigma(y-v) - \sigma$.

\item Set $\; k  \ \leftarrow \ K(z-w) + \lambda$.

\item Set
\[
\hspace{-0.1cm}
\left(
\hspace{-0.05cm}
\begin{array}{c}
y^+ \\[0.1cm]
z^+ \\[0.1cm]
\lambda^+
\end{array}
\hspace{-0.05cm}
\right)  \leftarrow \left(
\hspace{-0.05cm}
\begin{array}{rrr}
H & 0 & E^\tr \\[0.1cm]
0& K & -\mathbb{1} \\[0.1cm]
E & -\mathbb{1} & 0 
\end{array}
\hspace{-0.05cm}
\right)^{-1} \hspace{-0.1cm} \left(
\begin{array}{c}
Hv - h \\[0.1cm]
Kw - k \\[0.1cm]
0
\end{array}
\right)
\]

\item Set $\; y \ \leftarrow \ y^+$, $\; z \ \leftarrow \ z^+$, and $\; \lambda \ \leftarrow \ \frac{3}{4} \lambda^+ + \frac{1}{4} k$.

\item If $\; \log_3(i) \in \mathbb N \;$, do the following:

\begin{enumerate}

\item Attempt to guess the active set based on the current iterate. If successful, perform a single active set step and return the optimal solution.

\item Set $r = \frac{11}{10} \cdot \| w-z\|_\infty \ $, set
$$t = \frac{1}{\max \{  \| w-z\|_\infty, \| \mathcal Q y + \sigma \|_\infty  \}} \ ,$$
and update $\; K \ \leftarrow \ \nabla_z^{2} \Phi(r,t,z)$.

\end{enumerate}

\end{enumerate}

\item[] \textbf{End}

\smallskip

\hrule

\smallskip

\item \textbf{Output:} Primal and dual solution, $y$ and $\lambda$.

\smallskip

\hrule
\hrule
\end{itemize}

\caption{\label{fig::alg1} Summary of a version of ALADIN that is tailored for solving sparse QPs. All sparse linear algebra operations can be distributed by exploiting the sparse block structures of the matrices $E$ and $\mathcal Q$.}

\end{figure}
Notice that the first step of this algorithm computes the matrix-vector product $E^\tr \lambda$. Steps~2) and~3) set $v$ and $w$ to the explicit optimal solution of~\eqref{eq::dco1} and~\eqref{eq::dco2}, respectively. As $\mathcal Q + \Sigma$ is constant, its decomposition can be pre-computed exploiting the sparse and diagonal block-structure of $\mathcal Q$. Moreover, Step~4) evaluates the primal and dual residuums of the current iterates terminating the loop as soon as a user-specified termination tolerance $\mathsf{TOL} > 0$ is reached. Next, Steps~5) and~6) compute the subgradients in~\eqref{eq::gradF} and~\eqref{eq::gradG} while Step~7) computes the primal and dual solutions of the equality constrained consensus QP~\eqref{eq::consensus}. Although the presented ALADIN variant convergences for any value $\theta \in [0,1]$, the value $\theta = \frac{3}{4}$, as used by Step 8), has been found to work well on benchmark problems.\footnote{The parameter $\theta$ has here been tuned by empirical testing with thousands of randomly generated large-scale QPs. Adjusting this parameter properly leads to approximately $50\%$ of run-time improvement on average.}

Last but least, since Step 9) is computationally expensive, it is only evaluated whenever the iteration index is an integer power of $3$. In detail, Step 9a) attempts to guess the active set based on the active set of the projection step. If this yields a solution to the QP, one can directly terminate~\cite{Stellato2020}. Moreover, Step 9b) implements the log-barrier based Hessian approximation heuristic from the previous section. We use the barrier relaxation $r = \frac{11}{10} \| w-z \|_\infty$, a value that is a bit larger than $\| w-z \|_\infty$, such that the log-barrier is well-defined at the current iterate. Moreover, the log-barrier parameter $t$ is set to the inverse of the maximum of the primal and dual residuum of the current iterate. Notice that updating $K$ is expensive in the sense that one needs to update the sparse Cholesky decomposition of the KKT matrix in Step~7). However, as we will show below, this log-barrier update heuristic leads to significant overall run-time improvements.

\begin{remark}
Although the above variant of ALADIN has many similarities with ADMM~\cite{Stellato2020,Boyd2011}, it is not equivalent to ADMM, not even if we would skip the Hessian updates. It is, however, possible to construct variants of ALADIN that are equivalent to ADMM~\cite{Houska2016}. Nevertheless, the unique feature of ALADIN compared to ADMM is that it offers a natural way of choosing the augmented Lagrangian weights $H$ and $K$. Namely, we can exploit the similarity of ALADIN to SQP and interior point methods, which motivates the above log-barrier based scaling heuristic. As we will see below, the above outlined ALADIN variant performs well for large-scale QPs---even without using a pre-conditioner.
\end{remark}

\section{Real-Time Variants}
\label{sec::RealTime}

This section reviews ideas from~\cite{Jiang2021}, which can be used to develop a real-time ALADIN solver for distributed MPC. The performance of the solver from Figure~\ref{fig::alg1} and the performance of its real-time variant are discussed in Section~\ref{sec::results}.

\subsection{Real-Time Parallel MPC}
A real-time variant of the ALADIN based QP solver from Figure~\ref{fig::alg1} for the MPC problem~\eqref{eq::MPC} can be obtained by 
\begin{enumerate}
\addtolength{\itemsep}{2pt}
\item implementing only a finite number $i_\mathrm{max}$ of ALADIN iterations per sampling time,
\item sending the approximately optimal input $\hat u_0 \approx u_0^\star(x_0)$ to the real process as soon as $i = i_\mathrm{max}$,
\item skipping Step 9) but still updating $K$ once during every real-time loop using the relaxed log-barrier, and,
\item warm-starting the solver by shifting the solution from the previous time step. Additionally, the initialization variables are scaled such that $\| (y,z,\lambda) \| \leq \gamma_0 \| x_0 \|_Q$ for a sufficiently large constant $\gamma_0 < \infty$. 
\end{enumerate}
The above real-time ALADIN method has the property~\cite{Jiang2021} that there exist constants $\gamma_1 < \infty$ and $\kappa < 1$ with
\begin{align}
\label{eq::convergence}
\| x_0^+ - x_1^\star \| \leq \gamma_1  \kappa^{i_\mathrm{max}} \|  x_0 \|_Q \; .
\end{align}
Here, $\hat u_0 \approx u_0^\star(x_0)$ denotes the current ALADIN iterate for the first control input, $u_0^\star(x_0)$ the optimal input, \mbox{$x_0^+ = A x_0 + B \hat u_0$} the next state of the closed-loop system and \mbox{$x_1^\star = A x_0 + B u_0^\star$} the optimal solution for the state at the next time instance. Inequality~\eqref{eq::convergence} holds because ALADIN converges linearly. Moreover, due to our warm-start, $\Vert x_1^\star \Vert_Q$ scales at most linearly with respect to $\Vert x_0 \Vert_Q$; that is,
\[
\Vert x_1^\star \Vert_Q \leq \gamma_2 \| x_0 \|_Q
\]
for a constant $\gamma_2 < \infty$. Since the infinite horizon cost $J_\infty$ is piecewise quadratic, there exist constants \mbox{$\gamma_3,\gamma_4 < \infty$} with
\[
J_\infty(x)-J_\infty(x') \leq \gamma_3 \| x \|_Q \| x-x' \|_Q + \gamma_4 \| x-x' \|_Q^2
\]
for all $x,x'$ in the feasible domain of $J_\infty$. Next, by starting with the Bellman equation~\eqref{eq::Bellman} and substituting the three latter inequalities one finds that
\[
J_\infty(x_0^+) \leq J_\infty(x_0) - \| x_0 \|_Q^2 + [\gamma_2 \gamma_3   + \gamma_1 \gamma_4 ] \gamma_1 \kappa^{i_\mathrm{max}} \| x_0 \|_Q^2.
\]
This is a Lyapunov descent condition as long as
\[
i_\mathrm{max} > \frac{\log( \gamma_1 [\gamma_2 \gamma_3   + \gamma_1 \gamma_4 ] )}{ \log(\kappa^{-1}) } \; .
\]
Thus, if $i_\mathrm{max}$ is sufficiently large, the above real-time ALADIN variant yields an asymptotically stable feedback law as long as the iterates do not leave the domain of $J_\infty$~\cite{Jiang2021}. 

\begin{remark}
If the iterates of the real-time ALADIN variant leave the domain of $J_\infty$, the above estimates are wrong. As such, the above Lyapunov descent condition only ensures asymptotic stability under the assumption that recursive feasibility holds. If one is interested in a rigorous guarantee of recursive feasibility of real-time MPC in the presence of state constraints one needs to use methods from the field of rigid robust MPC~\cite{Langson2004} in order to pre-compute robustness margins for all state-constraints such that recursive feasibility holds in the presence of sufficiently small numerical errors.
\end{remark}

\begin{figure*}[t]
\begin{center}
\includegraphics[scale=0.45]{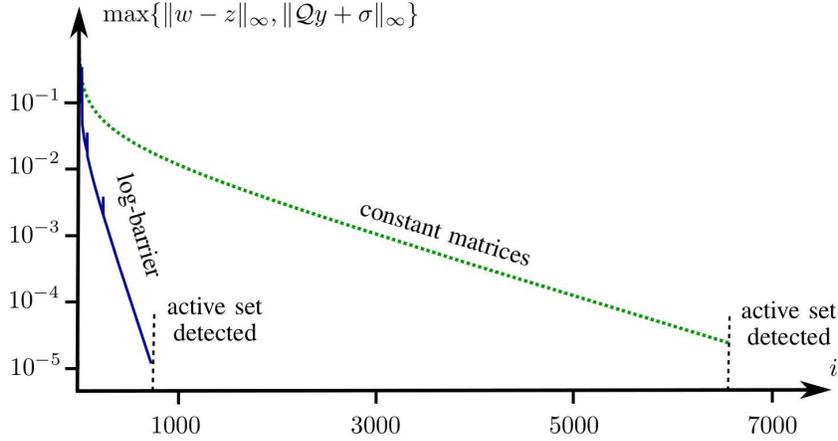}
\end{center}
\caption{\label{ref::convergence} Convergence of the ALADIN based sparse QP solver from Figure~\ref{fig::alg1} for the MPC problem~\eqref{eq::MPC} with and without log-barrier Hessian updates. The dotted green line shows the maximum of the primal and dual residuums versus the iteration index~$i$ of the main ALADIN loop for the case that the Hessian matrices are only updated once. In this case, the optimal active set is found after $6561$ iterations. The convergence profile can be compared to the solid blue line, which shows the ALADIN iteration for the case that the Hessian matrix $K$ is updated by using the log-barrier approach. The ``spikes'' in the convergence profile around the iterations $27$, $81$ and $243$ are caused by the Hessian updates, which can eventually lead to a drop of residual accuracy for a couple of iterations before paying out on the long run. In this example, the ALADIN based algorithm with log-barrier updates detects the optimal active set during iteration $729$.}
\end{figure*}

\subsection{Numerical Performance}
\label{sec::results}
In order to illustrate numerical performance, we implement the algorithm from Figure~\ref{fig::alg1} for the tutorial case study from Section~\ref{sec::tutorial}. We use $\mathsf{n} = 50$ wagons. This leads to an MPC problem with $100$ states and $50$ controls. We additionally set the prediction horizon to $N=100$. The corresponding QP is sparse: it has $15000$ optimization variables, $24900$ constraints, as well as $124202$ total non-zero entries in the QP data matrices and vectors. Figure~\ref{ref::convergence} shows the maximum of the primal and dual residuum versus the iteration index for the current state measurement $x_0 = 2 \cdot \mathbf{1}$. The log-barrier Hessian updates improve the convergence rate and overall run-time of the algorithm approximately by a factor~$8$.

In this case study, we implemented the presented ALADIN method in approximately $200$ lines of prototype Julia code finding that---for randomly chosen initial values $x_0$---this implementation needs on average $0.9$ seconds to solve the complete QP (without warm starts). On the same computer, \texttt{OSQP} solves the same QP in $0.8$ seconds on average while \texttt{GUROBI} needs more than $1$ second on average. The same trend in terms of run-time is confirmed by running these solvers on randomly generated sparse QPs. We do not elaborate more on this run-time result, because our goal here is merely to show that a simple implementation of ALADIN can achieve run-times that have the same order of magnitude as the run-times of existing sparse QP solvers.

\begin{figure}
\begin{center}
\includegraphics[scale=0.2]{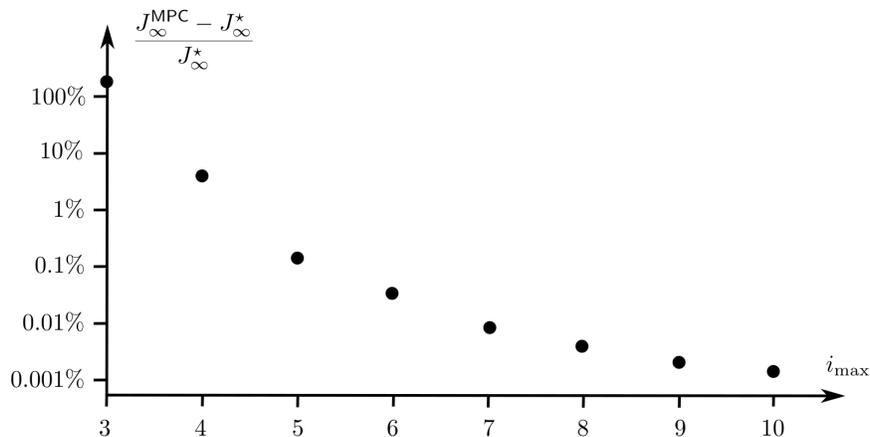}
\end{center}
\caption{\label{fig::performance} Relative loss of control performance of closed-loop real-time ALADIN compared to the optimal infinite horizon performance $J_\infty^\star = J_\infty(x_0)$ for the particular initial value $x_0 = 2 \cdot \mathbf{1}$ versus the maximum number $i_\mathrm{max}$ of ALADIN iterations per sampling time.}
\end{figure}

Finally, Figure~\ref{fig::performance} shows the loss of control performance,
\[
\frac{J_\infty^{\mathsf{MPC}} - J_\infty^\star}{J_\infty^\star} \quad \text{with} \quad J_\infty^{\mathsf{MPC}} = \sum_{k=0}^\infty \ell \left( x_k^\mathsf{MPC},\tilde \mu \left(x_k^\mathsf{MPC} \right) \right),
\]
where $\tilde \mu$ denotes the real-time ALADIN based MPC feedback law, $x_k^\mathsf{MPC}$ the associated approximately optimal closed loop trajectory and $J_\infty^\star = J_\infty(x_0)$ the optimal infinite horizon performance, both for the initial state $x_0 = 2 \cdot \mathbf{1}$. For $i_\mathrm{max} < 3$, the real-time ALADIN iteration happens to lead to an unstable closed-loop feedback law. However, for instance, for $i_\mathrm{max} = 5$ the relative loss of performance is smaller than $0.1 \%$. The run-time of real-time MPC is in this example less than $10 \, $ milliseconds---a run-time improvement of a factor $100$ compared to exact MPC.

\section{Conclusions}
\label{sec::conclusion}

This paper has presented a tutorial on how to implement a relatively simple ALADIN variant for solving sparse large-scale QPs as arising in the context of distributed MPC. It has been explained that this methods combines ideas from the field of sequential quadratic programming, interior point methods and augmented Lagrangian methods. For instance, in the proposed implementation, a relaxed log-barrier heuristic has been introduced in order to update certain Hessian matrices, which improves the convergence rate of ALADIN almost by an order of magnitude. 

The numerical results of this paper are relevant for the future of large-scale optimization and MPC solver development. This is because the presented ALADIN scheme can also be applied to solve more general convex as well as non-convex optimization problems~\cite{Houska2016,Engelmann2022}. The fact that the presented ``simple'' variant of this method can directly be used to implement a surprisingly competitive sparse QP solver within just a few lines of code can be interpreted as a promising indicator that ALADIN has enormous potential to be among the most competitive algorithms for large scale optimization and MPC. Besides, as pointed out in this tutorial, too, the algorithmic framework of ALADIN offers a unified perspective on augmented Lagrangian, SQP, and interior point methods. This perspective might help to proceed in a systematic way when synthesizing future large-scale non-convex optimization algorithms and software.

\newpage

\bibliographystyle{plain}
\bibliography{ALADIN}

\begin{thebibliography}{10}

\bibitem{Boyd2011}
S.~Boyd, N.~Parikh, E.~Chu, B.~Peleato, and J.~Eckstein.
\newblock Distributed optimization and statistical learning via the alternating
  direction method of multipliers.
\newblock {\em Found.~\& Trends in Machine Learning}, 3:1--122, 2011.

\bibitem{Engelmann2022}
A.~Engelmann, Y.~Jiang, H.~Benner, R.~Ou, B.~Houska, and T.~Faulwasser.
\newblock {ALADIN}-$\alpha$---an open-source {MATLAB} toolbox for distributed
  non-convex optimization.
\newblock {\em Optimal Control Applications \& Methods}, 43:4--22, 2022.

\bibitem{Ferreau2014}
H.J. Ferreau, C.~Kirches, A.~Potschka, H.G. Bock, and M.~Diehl.
\newblock qpoases: a parametric active-set algorithm for quadratic programming.
\newblock {\em Mathematical Programming Computation}, 6(4):327--363, 2014.

\bibitem{Frank1956}
M.~Frank and P.~Wolfe.
\newblock An algorithm for quadratic programming.
\newblock {\em Naval Res.~Log.~Q.}, 3:95--110, 1956.

\bibitem{Gertz2003}
E.M. Gertz and S.J. Wright.
\newblock Object-oriented software for quadratic programming.
\newblock {\em ACM Trans.~on Math.~Software}, 29(1):58--81, 2003.

\bibitem{Hamdi2011}
A.~Hamdi and S.K. Mishra.
\newblock Decomposition methods based on augmented {L}agrangian: a survey.
\newblock In {\em Topics in Nonconvex Optimization. Mishra, S.K., Chapter 11},
  pages 175--204, 2011.

\bibitem{Houska2011}
B.~Houska, H.J. Ferreau, and M.~Diehl.
\newblock An auto-generated real-time iteration algorithm for nonlinear {MPC}
  in the microsecond range.
\newblock {\em Automatica}, 47:2279–2285, 2011.

\bibitem{Houska2016}
B.~Houska, J.~Frasch, and M.~Diehl.
\newblock An augmented {L}agrangian based algorithm for distributed non-convex
  optimization.
\newblock {\em {SIAM} Journal on Optimization}, 26(2):1101--1127, 2016.

\bibitem{Houska2021}
B.~Houska and Y.~Jiang.
\newblock Distributed optimization and control with {ALADIN}.
\newblock {\em Recent Advances in Model Predictive Control: Theory, Algorithms,
  and Applications}, pages 135--163, 2021.

\bibitem{Jiang2021}
Y.~Jiang, J.~Oravec, B.~Houska, and M.~Kvasnica.
\newblock Parallel {MPC} for linear systems with input constraints.
\newblock {\em IEEE Transactions on Automatic Control}, 66(7):3401--3408, 2021.

\bibitem{Langson2004}
W.~Langson, I.~Chryssochoos, S.~V. Rakovi{\'c}, and D.~Q. Mayne.
\newblock Robust model predictive control using tubes.
\newblock {\em Autom.}, 40(1):125--133, 2004.

\bibitem{Mattingley2012}
J.~Mattingley and S.~Boyd.
\newblock {CVXGEN}: a code generator for embedded convex optimization.
\newblock {\em Optimization in Engineering}, 13(1):1--27, 2012.

\bibitem{Mosek2022}
MOSEK.
\newblock The {MOSEK} optimization toolbox for {MATLAB}, 2022.
  (\url{http://www.mosek.com}).

\bibitem{Muller2017}
M.A. M\"uller and F.~Allg\"ower.
\newblock Economic and distributed model predictive control: Recent
  developments in optimization-based control.
\newblock {\em Journal of Control, Measurement, and System Integration},
  10(2):39--52, 2017.

\bibitem{Nesterov1994}
Y.~Nesterov and A.~Nemirovskii.
\newblock {\em Interior-Point Polynomial Algorithms in Convex Programming}.
\newblock SIAM, Philadelphia, 1994.

\bibitem{Nocedal2006}
J.~Nocedal and S.J. Wright.
\newblock {\em Numerical Optimization}.
\newblock Springer Series in Operations Research and Financial Engineering
  Springer, 2006.

\bibitem{Donoghue2016}
B.~O'Donoghue, E.~Chu, N.~Parikh, and S.~Boyd.
\newblock Conic optimization via operator splitting and homogeneous self-dual
  embedding.
\newblock {\em Journal of Optimization Theory and Applications},
  169(3):1042--1068, 2016.

\bibitem{Gurobi2022}
Gurobi Optimization.
\newblock Gurobi optimizer reference manual, 2022.
  (\url{http://www.gurobi.com}).

\bibitem{Qin2003}
S.J. Qin and T.A. Badgwell.
\newblock A survey of industrial model predictive control technology.
\newblock {\em Con.~Eng.~Practice}, 93(316):733--764, 2003.

\bibitem{Rawlings2017}
J.B. Rawlings, D.Q. Mayne, and M.M. Diehl.
\newblock {\em Model predictive control: Theory and design}.
\newblock Nob Hill Publishing, 2017.

\bibitem{Stellato2020}
B.~Stellato, G.~Banjac, P.~Goulart, A.~Bemporad, and S.~Boyd.
\newblock {OSQP}: an operator splitting solver for quadratic programs.
\newblock {\em Mathematical Programming Computation}, 12:637--672, 2020.

\bibitem{Wolfe1959}
P.~Wolfe.
\newblock The simplex method for quadratic programming.
\newblock {\em Econometrica}, 27(3):382--398, 1959.

\end{thebibliography}

\end{document}